# MODULARITY PREVENTS TAILS

KEITH A. KEARNES AND EMIL W. KISS

ABSTRACT. We establish a direct correspondence between two congruence properties for finite algebras. The first property is that minimal sets of type **i** omit tails. The second property is that congruence lattices omit pentagons of type **i**.

## 1. INTRODUCTION

Tame congruence theory is a framework developed to study the local polynomial structure of finite algebras. The creed of the tame congruence theorist is that it is always useful to localize to a small subset, and then to use combinatorial methods to understand how operations behave on that subset. These small subsets are called 'minimal sets', and even on these minimal sets tame congruence theory can only give partial information. Specifically, one finds in [2] that minimal sets are divided into two main parts, the 'body' and the 'tail'. On the body the theory gives almost complete information about the structure, while on the tail the theory gives almost no information. The body of a minimal set is always nonempty although the tail may be empty. When the tail is empty celebrations usually ensue because the theory can be used to its utmost.

In this paper we take a bad situation (minimal sets with tails), and then make it worse! We do this by showing that there must be other minimal sets nearby which have more complicated tails, tails which are so complicated that their existence can be detected by simply looking at the shape of a congruence lattice. But this is really looking at the situation backwards: the reason for proving the result we do is to show how restrictions on the labeled congruence lattices of algebras in a variety force the nonexistence of tails in minimal sets.

Let $\mathbf{A}$ be a finite algebra. The precise statement we prove about $\mathbf{A}$, and every type $\mathbf{i} \in \{\mathbf{1}, \mathbf{2}, \mathbf{3}, \mathbf{4}, \mathbf{5}\}$, is the following. *If, for all subalgebras $\mathbf{S} \leq \mathbf{A}^2$, the labeled congruence lattice of $\mathbf{S}$ has no pentagon whose critical quotient is of type $\mathbf{i}$, then the minimal sets of type $\mathbf{i}$ in $\mathbf{A}$ have empty tails.* In particular, if the congruence lattices of subalgebras of $\mathbf{A}^2$ are modular, then all minimal sets of $\mathbf{A}$ have empty tails. This

1991 *Mathematics Subject Classification.* 08A05, 08A30, 08B10.
Work supported by the Hungarian National Foundation for Scientific Research, grant no. 16432, and by the Fields Institute (Toronto, Canada).

is an improvement of Theorem 8.5 of [2], where the same conclusion is proved under the hypothesis that the variety generated by **A** is congruence modular.

In this paper sets will be denoted by capital letters in italics, as in $S, T, U, \ldots$, and algebras in bold face capitals, as in $\mathbf{A}, \mathbf{B}, \mathbf{C}, \ldots$. Elements of sets are usually denoted by lower case letters in italics, as in $x, a, b, c, \ldots$. A lower case letter in bold face denotes a sequence of some length, so $\mathbf{x}$ denotes something like $(x_1, x_2, \ldots, x_m)$. If **A** is an algebra, and we are investigating $\mathbf{A}^n$, then $\hat{a}$ will denote the constant sequence $(a, \ldots, a) \in A^n$. The set of all these constant sequences for $a \in A$ is called the diagonal of $A^n$. If $p(\mathbf{x})$ is a polynomial of **A**, then for any subalgebra $\mathbf{S} \leq \mathbf{A}^n$ containing the diagonal, the operation which is $p$ in each coordinate is a polynomial of **S** which will be denoted $\hat{p}$.

## 2. Modularity and Tails

If $\delta \prec \theta$ is a covering pair of congruences on **A**, then a *minimal set for the quotient* $\langle \delta, \theta \rangle$ is a set $U \subseteq A$ which is minimal under inclusion among all sets of the form $e(A)$, where $e$ is an idempotent unary polynomial of **A**, for which the restrictions $\delta|_U$ and $\theta|_U$ are different. A *trace* of $U$ is a $\theta|_U$-class which differs from a $\delta|_U$-class. The *body* $B$ of $U$ is the union of the traces, and the *tail* is $T = U - B$. In this section we show that if $T$ is nonempty, then the minimality of $U$ forces a failure of congruence modularity in some subalgebra of $\mathbf{A}^2$.

Minimal sets come in different types depending on the kind of polynomials of **A** which can be restricted to the minimal set. This type is indicated by assigning a type label to the minimal sets for a quotient $\langle \delta, \theta \rangle$, $\delta \prec \theta$, and to the quotient itself. It is explained in [2] how to do this. The type label is a number $\mathbf{i} \in \{\mathbf{1}, \mathbf{2}, \mathbf{3}, \mathbf{4}, \mathbf{5}\}$, and we write $\mathrm{typ}(\delta, \theta) = \mathbf{i}$ to indicate the label assigned to $\langle \delta, \theta \rangle$.

A *pentagon* $[\gamma, \delta, \theta]$ in a lattice is a five element sublattice generated by elements $\gamma, \delta, \theta$ satisfying $\delta < \theta$, $\gamma \vee \delta \geq \theta$ and $\gamma \wedge \theta \leq \delta$ (see Figure 1).

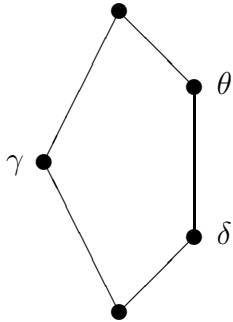

FIGURE 1: The pentagon $[\gamma, \delta, \theta]$

A lattice is *modular* if no sublattice is a pentagon. The *critical quotient* in the pentagon $[\gamma, \delta, \theta]$ is the quotient $\langle \delta, \theta \rangle$. Observe that if $[\gamma, \delta, \theta]$ is a pentagon with critical quotient $\langle \delta, \theta \rangle$ and $\delta \leq \delta' < \theta' \leq \theta$, then $[\gamma, \delta', \theta']$ is a pentagon with critical quotient $\langle \delta', \theta' \rangle$. Thus, as we will concern ourselves with pentagons in congruence lattices of finite algebras, we can always shrink the critical quotient to a covering quotient. This covering has a type label, as mentioned above, and if that type label is **i** we will call the pentagon a *pentagon of type* **i**.

The easy part of the connection between failures of congruence modularity and empty tails is stated in the following theorem.

**Theorem 2.1.** *Assume that* **A** *is a finite algebra and* $\mathrm{Con}(\mathbf{A})$ *has a pentagon* $[\gamma, \delta, \theta]$ *of type* **i**. *If* **i** $\neq$ **1** *and* $U$ *is a* $\langle \delta, \theta \rangle$*-minimal set, then* $U$ *has nonempty tail.*

*Proof.* In fact we will show that if $B$ and $T$ are the body and tail of $U$, then
$$\gamma \cap (B \times T) \neq \emptyset.$$
If instead $\gamma \cap (B \times T) = \emptyset$, then the body $B$ is a union of $\gamma|_U$-classes. From the definition of the body, $B$ is a union of $\theta|_U$-classes. Thus, if the theorem is false then $B$ is a union of $\gamma \vee \delta$-classes. If this is so, then restriction to the body is a lattice homomorphism from the pentagon $[\gamma, \delta, \theta]$ to a sublattice generated by $\{\gamma|_B, \delta|_B, \theta|_B\}$. This lattice homomorphism doesn't identify $\delta$ and $\theta$ since the definition of the body prevents $\delta|_B$ equaling $\theta|_B$. Thus $[\gamma|_B, \delta|_B, \theta|_B]$ is a pentagon in $\mathrm{Con}(\mathbf{A}|_B)$. Now if the type is **2**, **3** or **4** we have a contradiction, since $\mathrm{Con}(\mathbf{A}|_B)$ is congruence modular. In types **3**, **4** or **5** we get a different contradiction: $B$ is a $\theta|_B$-class, so $B$ cannot support a pentagon of the form $[\gamma|_B, \delta|_B, \theta|_B]$. This finishes the proof. □

The next lemma is a necessary tool for the proof of our main result.

**Lemma 2.2.** *Let* **C** *be a finite* $\langle \alpha, \beta \rangle$*-minimal algebra whose body and tail are $B$ and $T$ respectively. Assume that $t \in T$ and $p(\mathbf{x})$ is an $m$-ary polynomial of* **C** *which satisfies that* $p(t, t, \ldots, t) \in B$. *Then* $p(\beta|_B, \beta|_B, \ldots, \beta|_B) \subseteq \alpha$.

*Proof.* There is no loss of generality by assuming $\alpha = 0$, so we do. Suppose that the statement of the lemma fails for some $p$ and $t$. Our first aim is to construct a failure of the lemma with a binary polynomial.

**Claim 2.3.** *There is an element $0 \in B$ and a binary polynomial $q(x, y)$ such that*
(1) $q(t, t) \in B$; *and*
(2) $q(x, 0) = x$ *holds on $C$.*

*Proof.* We show how to modify $p$ to get $q$. First, since $p(\beta|_B, \beta|_B, \ldots, \beta|_B) \not\subseteq \alpha$, there exist $\mathbf{a}, \mathbf{b} \in B^m$ such that $(a_i, b_i) \in \beta$ for all $i$ and $p(\mathbf{a}) \not\equiv_\alpha p(\mathbf{b})$. If we define
$$\mathbf{c}_i := (b_1, b_2, \ldots, b_{i-1}, a_i, \ldots, a_m),$$
then since $p(\mathbf{c}_1) = p(\mathbf{a}) \not\equiv_\alpha p(\mathbf{b}) = p(\mathbf{c}_{m+1})$ we must have $p(\mathbf{c}_i) \not\equiv_\alpha p(\mathbf{c}_{i+1})$ for some $i$. Replacing **a** and **b** with $\mathbf{c}_i$ and $\mathbf{c}_{i+1}$ we see that we may assume that $a_j = b_j$ for all but one value of $j$. We will assume that $\mathbf{a}$ and $\mathbf{b}$ differ in the first component only, so $p(\mathbf{a}) = p(a_1, a_2, \ldots, a_m) \not\equiv_\alpha p(b_1, a_2, \ldots, a_m) = p(\mathbf{b})$.

We have shown that there is some $\mathbf{u} \in B^{m-1}$, namely $\mathbf{u} = (a_2, \ldots, a_m)$, such that $p(\beta|_B, \mathbf{u}) \not\subseteq \alpha$. This proves that $p(x, \mathbf{u})$ is a permutation of $C$. To produce from $p$ a binary polynomial $q$ which satisfies the conditions of the claim, we will separate into two cases which depend on whether or not $\langle \alpha, \beta \rangle$ is abelian or nonabelian.

If this quotient is abelian, which means that the type is **1** or **2**, we know by Theorem 3.4 of [8] that the body $B$ is contained in a single class of the twin congruence of $\mathbf{C}$. This means that, since $p(x, \mathbf{u})$ is a permutation of $C$ and $\mathbf{u} \in B^{m-1}$, therefore $p(x, \mathbf{v})$ is a permutation for any $\mathbf{v} \in B^{m-1}$. So, choose any $0 \in B$ and let $p'(x, y) = p(x, y, \ldots, y)$. Then $p'(t, t) = p(t, t, \ldots, t) \in B$ and $p'(x, 0)$ is a permutation of $C$.

To get a similar $p'$ in types **3**, **4** and **5**, let $+$ be a semilattice polynomial on $B$ for which there is a $0 \in B$ satisfying $0 + x = x + 0 = x$ on $C$. There exist such $0$ and $+$ by Lemmas 4.15 and 4.17 of [2]. Now let $p'(x, y) = p(x, y + u_1, \ldots, y + u_m)$. This polynomial satisfies that $p'(x, 0) = p(x, \mathbf{u})$ is a permutation of $C$ and, since $B^2 \subseteq \beta$,

$$p'(t, t) = p(t, t + u_1, \ldots, t + u_m) \equiv_\beta p(t, t + 0, \ldots, t + 0) = p(t, t, \ldots, t) \in B.$$

For all five types the polynomial $p'$ has all the desired properties of $q$, except that we have to strengthen the fact that $p'(x, 0)$ is a permutation to $q(x, 0) = x$. So, if $\pi(x) = p'(x, 0)$, then some power $\pi^{-1}$ acts like the inverse of $\pi$ on $C$. Define $q(x, y) := \pi^{-1}(p'(x, y))$. This polynomial clearly satisfies (2). Since permutations of $C$ restrict to permutations of $B$ and of $T$, we have $\pi(B) = B$ and so the fact that $p'$ satisfied (1) implies that $q$ does too. This completes the proof of Claim 2.3. ∎

During the rest of the argument we will let $N$ be the $\langle \alpha, \beta \rangle$-trace of $\mathbf{C}$ which contains the element 0 from Claim 2.3.

Assume first that the type of $\langle \alpha, \beta \rangle$ is **1**. Since $q(x, 0) = x$ on $C$, the trace $N$ is closed under $q$ and $q$ depends on its first variable on $N$. Hence, $q$ does not depend on its other variables on $N$, therefore we have

$$h(n) := q(n, n) = q(n, 0) = n$$

for $n \in N$. Since $\mathbf{C}$ is an $\langle \alpha, \beta \rangle$-minimal algebra, it follows that $h$ is a permutation of $C$. However, restriction (1) forces

$h(t) = q(t, t) \in B$, and no permutation of $C$ can map a tail element into the body. This contradiction completes the type **1** argument.

In the other types we shall use a similar argument, but our unary polynomial is

$$H(x) := q(x, q(t, q(t, x))).$$

We shall prove that $H(t) \in B$, so $H$ cannot be a permutation, and will get the desired contradiction by calculating $H$ on the trace $N$, and show that is one-to-one on $N$.

Let $b = q(t, t) \in B$. Then $H(t) = q(t, q(t, b))$. In all types different from **1** there exists a congruence $\rho$ of $\mathbf{C}$ such that $B$ is a $\rho$-class: if the type is **2** we can take $\rho$ to be the twin congruence, while in the nonabelian types $\rho = \beta$ works, since $B = N$. Hence $q(t, 0) = t$ implies that

$$H(t) = q(t, q(t, b)) \equiv_\rho q(t, q(t, 0)) = q(t, t) \in B.$$

Thus $H(t) \in B$ indeed.

Now let us compute $H(n)$ for $n \in N$. First notice that $q(t, n) \equiv_\beta q(t, 0) = t$. As $t$ is in the tail, $t/\beta = t/\alpha = \{t\}$, and so we get that $q(t, n) = t$. Therefore $H(n) = q(n, q(t, t)) = q(n, b)$. We have to show that this, as a function of $n$, is one-to-one on $N$.

This is clear if the type of $\mathbf{C}$ is $\mathbf{2}$, since $q(n, 0) = n$, and $0 \in B$ is related to $b \in B$ by the twin congruence. Thus we can assume that the type is nonabelian. In that case, we have that $N = \{0, 1\}$. From $q(x, 0) = x$ we see that in the case when $b = 0$ we are done. So we can assume that $b = 1$, and we can finish the proof by showing that $q(0, 1) \neq q(1, 1)$. In fact, we shall prove that $q(1, 1) = 0$ and $q(0, 1) = 1$.

Consider again the polynomial $h(x) := q(x, x)$. Since $h(t) = b \in B$, this is not a permutation of $C$. Therefore, we have

$$q(1, 1) = h(1) = h(0) = q(0, 0) = 0.$$

Next we show that $q(0, 1) = 1$. From $q(0, 0) = 0$ we see that $q(N, N) \subseteq N$, so it is sufficient to prove that $q(0, 1) \neq 0$. Suppose that in fact $q(0, 1) = 0$, and consider the polynomial $k(x) := q(x, q(x, 1))$. We have

$$k(t) = q(t, q(t, 1)) \equiv_\beta q(t, q(t, 0)) = q(t, t) \in B.$$

This shows that $k(x)$ is not a permutation. Therefore

$$0 = q(0, 0) = q(0, q(0, 1)) = k(0) = k(1) = q(1, q(1, 1)) = q(1, 0) = 1.$$

This contradiction finishes the proof of Lemma 2.2. □

**Theorem 2.4.** *Assume that $\alpha \prec \beta$ in $\mathrm{Con}(\mathbf{A})$ and that $U$ is an $\langle \alpha, \beta \rangle$-minimal set of type $\mathbf{i}$ which has nonempty tail. Then there is a subalgebra $\mathbf{S} \leq \mathbf{A}^2$ whose congruence lattice contains a pentagon of type $\mathbf{i}$.*

*Proof.* Denote by $B$ and $T$ the body and tail of $U$ respectively. Let $t$ be an element of $T$ and choose $(0, 1) \in \beta|_B - \alpha|_B$. Define $\mathbf{S}$ to be the subalgebra of $\mathbf{A}^2$ generated by the set

$$\Delta \cup \{(0, 1), (0, t), (1, t)\},$$

where $\Delta$ denotes the diagonal of $\mathbf{A}^2$. Let $\gamma = \alpha \times 1_\mathbf{A}$ denote the congruence on $\mathbf{S}$ which relates pairs $(a, b)$ and $(c, d)$ if and only if $(a, c) \in \alpha$. Let $\delta$ be the congruence on $\mathbf{S}$ generated by

$$(\alpha \times \alpha) \cup \{\langle (0, t), (1, t) \rangle\}.$$

Let $\theta$ be the congruence generated by $\delta$ and the pair $\langle (0, 1), (1, 1) \rangle$. Clearly, we have $\alpha \times \alpha \leq \delta \leq \theta \leq 1_\mathbf{A} \times \alpha$. Here is a picture of this situation.

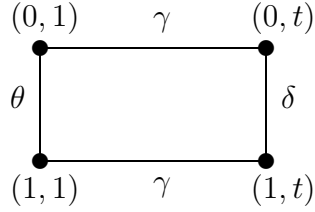

**Claim 2.5.** $[\gamma, \delta, \theta]$ *is a pentagon in* $\mathrm{Con}(\mathbf{S})$.

*Proof.* Since pairs in $\gamma$ are $\alpha$-related in the first coordinate and pairs in $\theta$ are $\alpha$-related in the second, we get that $\gamma \wedge \theta \leq \alpha \times \alpha \leq \delta$. Since
$$(0,1) \equiv_\gamma (0,t) \equiv_\delta (1,t) \equiv_\gamma (1,1)$$
we get that $\gamma \vee \delta \geq \theta$. So we only need to show that $\delta < \theta$ to finish the proof of the claim. Equivalently, we must show that $\langle (0,1), (1,1) \rangle \notin \delta$.

Let $e$ be an idempotent unary polynomial of $\mathbf{A}$ for which $e(A) = U$. Let $\hat{e}$ denote the polynomial of $\mathbf{S}$ which is $e$ acting coordinatewise. Let $V = \hat{e}(S) = U^2 \cap S$. Since $(0,1), (1,1) \in V$ and $V$ is the image of an idempotent unary polynomial, it follows that if $\langle (0,1), (1,1) \rangle \in \delta$ then there is a Maltsev chain lying entirely in $V$ which connects $(1,1)$ to $(0,1)$ by polynomial images of generating pairs for $\delta$. To prove that there is no such chain it will suffice to prove that if $\langle g, h \rangle$ is a generating pair for $\delta$ and $F$ is a unary polynomial of $\mathbf{S}$, then
$$\hat{e}F(g) \equiv_{\alpha \times \alpha} (1,1) \iff \hat{e}F(h) \equiv_{\alpha \times \alpha} (1,1).$$
For then we will have that the $\delta|_V$-class containing $(1,1)$ is the same as the $(\alpha \times \alpha)|_V$-class containing $(0,1)$. The latter does not contain $(0,1)$, so the former cannot.

Let us assume that the displayed bi-implication fails to hold. Then there is a unary polynomial $F$ of $\mathbf{S}$ and a generating pair $\langle g, h \rangle$ for $\delta$ such that $\hat{e}F(g) \equiv_{\alpha \times \alpha} (1,1)$ while $\hat{e}F(h) \not\equiv_{\alpha \times \alpha} (1,1)$. Necessarily $\langle g, h \rangle \notin (\alpha \times \alpha)$, and therefore we must have $\langle g, h \rangle = \langle (0,t), (1,t) \rangle$ or $\langle (1,t), (0,t) \rangle$. The argument in both cases is similar, so we assume that $g = (1,t)$ and $h = (0,t)$. The polynomial $\hat{e}F(X)$ may be expressed as
$$\hat{e}F(X) = \hat{r}(X, (0,t), (1,t), (0,1))$$
for some 4-ary polynomial $r$ of $\mathbf{A}$. Since $\hat{e}F$ is prefixed by $\hat{e}$, we get that $r(A^4) \subseteq U$. Thus, $r$ is a polynomial of the induced algebra $\mathbf{C} = \mathbf{A}|_U$. Since $\hat{e}F(0,t) \equiv_{\alpha \times \alpha} (1,1)$ and $\hat{e}F(1,t) \not\equiv_{\alpha \times \alpha} (1,1)$, we get that
- $r(1,0,1,0) \equiv_\alpha 1$;
- $r(0,0,1,0) \not\equiv_\alpha 1$; and
- $r(t,t,t,1) \equiv_\alpha 1$.

We define $p(x,y,z) = r(x,y,z,0)$. This is a polynomial of $\mathbf{C}$, and we have
$$p(t,t,t) = r(t,t,t,0) \equiv_\beta r(t,t,t,1) \equiv_\beta 1 \in B.$$

However, $(0, 1) \in \beta|_B$ and
$$p(0, 0, 1) = r(0, 0, 1, 0) \not\equiv_\alpha r(1, 0, 1, 0) = p(1, 0, 1).$$
This shows that $p(\beta|_B, \beta|_B, \beta|_B) \not\subseteq \alpha$. Lemma 2.2 proves that there is no such $p$, therefore we have $\delta < \theta$ and the claim is proven. ∎

We now show that it is possible to shrink the interval between $\delta$ and $\theta$ to a prime quotient whose type is **i**.

**Claim 2.6.** $\gamma \prec (\gamma \vee \theta)$ and this covering is of type **i**.

*Proof.* Clearly, $\mathbf{S} \leq \mathbf{A}^2$ is a subdirect representation and $\gamma$ is the coimage of $\alpha$ under the first coordinate projection. Both $\gamma$ and $\theta$ are contained in the coimage of $\beta$ under the first coordinate projection (which is equivalent to saying that all $\langle(a, b), (c, d)\rangle \in \gamma \cup \theta$ have the property that $(a, c) \in \beta$). Therefore, since $\gamma < (\gamma \vee \theta)$, we must have that $(\gamma \vee \theta)$ is exactly the coimage of $\beta$, so $\gamma \prec (\gamma \vee \theta)$ and this covering is of type **i**. This labeling is shown in Figure 2. ∎

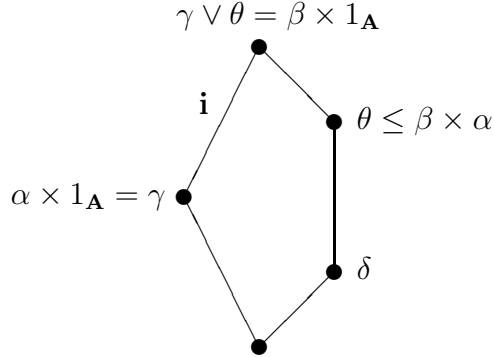

FIGURE 2: $\gamma \prec (\gamma \vee \theta)$ has the same type as $\alpha \prec \beta$.

Claim 2.5 produces a pentagon, but we need to shrink the critical quotient to a covering in order for it to qualify as a 'pentagon of type **i**' for some **i**. Therefore, replace $\delta$ by any $\delta'$ such that $\delta \leq \delta' \prec \theta$. This way we get a new pentagon $[\gamma, \delta', \theta]$. What we need to show now is that the type of $\langle \delta', \theta \rangle$ is the same as the type of $\langle \alpha, \beta \rangle$. For this task, note that since $\theta$ is generated by $\delta \cup \{\langle(0, 1), (1, 1)\rangle\}$ we have $\langle(0, 1), (1, 1)\rangle \in \theta - \delta'$. We will use this fact to determine the type of $\langle \delta', \theta \rangle$.

If **i** = **1**, then the quotient $\langle \gamma, \gamma \vee \theta \rangle$ is strongly abelian, which implies that $\theta$ is strongly abelian over $\gamma \wedge \theta$. Every prime quotient between $\gamma \wedge \theta$ and $\theta$ must have type **1**, so $\langle \delta', \theta \rangle$ is of type **1**.

If **i** = **2**, then by the same argument as above we can conclude that $\theta$ is abelian over $\gamma \wedge \theta$, so the type of $\langle \delta', \theta \rangle$ is either **1** or **2**. However, the pair $\langle(0, 1), (1, 1)\rangle$ is a **1**-snag, as is witnessed by the polynomial $b(X, Y) := \hat{d}(X, (1, 1), Y)$ where $d(x, y, z)$ is a pseudo-Maltsev polynomial for $U$. This shows that the type cannot be **1**, so it must be **2**.

Now we handle the case when **i** is **3**, **4** or **5**. Recall that $V = \hat{e}(S) = U^2 \cap S$. Since $\langle (0,1), (1,1) \rangle \in \theta|_V - \delta'$ we have $\delta'|_V < \theta|_V$. This implies that $V$ contains a $\langle \delta', \theta \rangle$-minimal set $W$; moreover there is a unary polynomial $P(X)$ of **S** such that $\hat{e}P(S) \subseteq W \subseteq V$ and $\hat{e}P((0,1)) \not\equiv_{\delta'} \hat{e}P((1,1))$. Expressing $\hat{e}P(X)$ as $\hat{p}(X, (0,t), (1,t), (0,1))$, we find that for $X = (x, y)$ the polynomial $\hat{e}P(X)$ is of the form $\langle p(x, \mathbf{u}), p(y, \mathbf{v}) \rangle$ where $p(A^4) \subseteq U$. If $p(0, \mathbf{u}) \equiv_\alpha p(1, \mathbf{u})$, then we would have

$$\langle \hat{e}P((0,1)), \hat{e}P((1,1)) \rangle \in (\alpha \times \alpha) \leq \delta \leq \delta',$$

which is false. Therefore $p(\beta|_U, \mathbf{u}) \not\subseteq \alpha|_U$, and so $p(x, \mathbf{u})$ is a permutation of $U$. Let $\pi(x)$ denote an iterate of $p(x, \mathbf{u})$ which acts like the inverse of this permutation on $U$. Then $\widehat{\pi}(X)$ is a permutation of $V$, so $\widehat{\pi}\hat{e}P(X)$ is a polynomial which maps the pair $\langle (0,1), (1,1) \rangle$ into a $\langle \delta', \theta \rangle$-minimal set in $V$, and the pair of elements it is mapped to has the form

$$\langle \widehat{\pi}\hat{e}P((0,1)), \widehat{\pi}\hat{e}P((1,1)) \rangle = \langle (0, w), (1, w) \rangle$$

where $w \in U$.

We started with the assumption that the type of **i** was **3**, **4** or **5** and that $(0,1)$ was an $\langle \alpha, \beta \rangle$-subtrace contained in the minimal set $U$. At this point we know that there is some $w \in U$ such that $\langle (0, w), (1, w) \rangle$ is a $\langle \delta', \theta \rangle$-subtrace. This is enough to show that the type of $\langle \alpha, \beta \rangle$ is the same as the type of $\langle \delta', \theta \rangle$. We explain now why this is so. Let **J** denote the algebra whose universe is $\{0, 1\}$ and whose operations are the idempotent (polynomial) operations of $\mathbf{A}|_{\{0,1\}}$. Let **K** denote the algebra whose universe is $\{(0, w), (1, w)\}$ and whose operations are the idempotent operations of $\mathbf{S}|_{\{(0,w),(1,w)\}}$. The bijection $x \mapsto (x, w)$ is a weak isomorphism from **J** to **K**, that is, it induces a bijection between the clones of these algebras. The inverse bijection is just projection onto the first component. To verify this claim, note that if $q(\mathbf{x})$ is an idempotent operation of $\mathbf{J} = \mathbf{A}|_{\{0,1\}}$, then $q$ is the restriction of a polynomial of **A** which is idempotent on $U$, so $\hat{q}(\mathbf{X})$ is a corresponding idempotent operation of $\mathbf{K} = \mathbf{S}|_{\{(0,w),(1,w)\}}$. Conversely, given any idempotent operation of **K**, one gets a corresponding idempotent operation on **J** by just projecting to the first component. These correspondences between the clones of **J** and **K** can be checked to be composition preserving and to be inverses of each other.

However, in the nonabelian cases the type of a quotient is determined by the idempotent polynomial operations of any subtrace, because join, meet, and ternary addition modulo 2 are idempotent operations on $\{0, 1\}$. Therefore, if **i** is **3**, **4** or **5**, then $\{0, 1\}$ supports the idempotent polynomial operations of a Boolean algebra, lattice or semilattice respectively. We have shown that $\{(0, w), (1, w)\}$ supports the same type of idempotent operations, so $\langle \delta', \theta \rangle$ is nonabelian and of the same type. This finishes the proof. □

As mentioned in the introduction, Theorem 2.4 implies that if all subalgebras of $\mathbf{A}^2$ are congruence modular, then all minimal sets of $\mathbf{A}$ have empty tail. With a little argument this can be refined.

**Corollary 2.7.** *Let $\mathbf{A}$ be a finite algebra. Then* $(1) \implies (2) \implies (3)$.
 (1) *All subalgebras of $\mathbf{A}^3$ have modular congruence lattices.*
 (2) *All minimal sets of $\mathbf{A}$ are of types $\mathbf{2}, \mathbf{3}$ or $\mathbf{4}$ and have empty tails.*
 (3) *$\mathbf{A}$ has a modular congruence lattice.*
*Moreover* $(a) \implies (b) \implies (c)$.
 (a) *All subalgebras of $\mathbf{A}^2$ have distributive congruence lattices.*
 (b) *All minimal sets of $\mathbf{A}$ are of types $\mathbf{3}$ or $\mathbf{4}$ and have empty tails.*
 (c) *$\mathbf{A}$ has a distributive congruence lattice.*

*Proof.* The previous theorem explains why (1) and (a) imply the empty tails part of (2) and (b). Theorem 5.27 of [2] shows that both (1) and (a) imply that $\mathbf{5}$ does not appear in the congruence lattice of $\mathbf{A}$. A similar argument works for type $\mathbf{1}$. Note though that one must work in a subalgebra of $\mathbf{A}^3$ to show that congruence modularity forces type $\mathbf{1}$ to be omitted, but it suffices to work in a subalgebra of $\mathbf{A}^2$ if one uses congruence distributivity to force $\mathbf{1}$ to be omitted. The proof of Theorem 8.7 of [2] shows that $(2) \implies (3)$, while a similar argument shows that $(b) \implies (c)$. □

We were motivated to consider the question of whether tails of type $\mathbf{i}$ force the existence of pentagons of type $\mathbf{i}$ for the following two reasons. The first reason is that the pentagon shape is connected with the centralizer relation by the fact that if $[\gamma, \delta, \theta]$ is a pentagon in a congruence lattice, then $C(\theta, \gamma; \delta)$ holds. Therefore, if a type $\mathbf{i}$ quotient appears in a variety, and its minimal sets have tails, then one can find such an example where the quotient is centralized by a congruence which connects the body to the tail. This is an (admittedly technical but) important fact. The second reason we considered this question is that it often happens that interesting algebraic conditions can only occur in varieties with restrictions on type labelings. We wanted a way to directly translate this information about restricted type labelings into properties of the minimal sets, these being the smallest sets which support good local approximation. The following sequence of corollaries explain how this is done.

In [1] it is shown that the congruence lattice of a finite algebra in a congruence semimodular variety has no pentagons of type $\mathbf{2}, \mathbf{3}$ or $\mathbf{4}$. This has the following consequence, which is new and answers a question posed in [1].

**Corollary 2.8.** *If $\mathcal{V}$ is a locally finite, congruence semimodular variety, then the minimal sets of types $\mathbf{2}$, $\mathbf{3}$ or $\mathbf{4}$ have empty tails.*

It is shown in [3] that any locally finite, residually small variety which omits type $\mathbf{5}$ has the property that there is no pentagon of type $\mathbf{2}, \mathbf{3}, \mathbf{4}$ (or $\mathbf{5}$). The following corollary to Theorem 2.4 has never been published, although it seems to have been known to some. (For example, results in [4] or [7] are older than this paper and they contain sufficient information to deduce the next corollary.)

**Corollary 2.9.** *If $\mathcal{V}$ is a locally finite, residually small variety which omits type* **5**, *then all minimal sets of types* **2**, **3** *or* **4** *have empty tails.*

Using the construction of subdirectly irreducible algebras in [6] and Claim 2.5 one can get more information on how the operations behave "on the border" between the body and the tail on minimal sets in residually small varieties, even if type **5** is present. For example, in a minimal set of type **2** there is no subset $\{b, t\}$ with $b$ in the body and $t$ in the tail which supports a meet operation with neutral element $t$.

In [5] it is shown that a locally finite variety has a difference term if and only if no finite member has a congruence lattice with a pentagon of type **1** or **2**. The following result is now an immediate corollary of this statement. It was also proved in [5] using different arguments.

**Corollary 2.10.** *If $\mathcal{V}$ is a locally finite variety, then $\mathcal{V}$ has a difference term if and only if it omits type* **1** *and the type* **2** *minimal sets have empty tails.*

## References

[1] P. Agliano and K. A. Kearnes, *Congruence semimodular varieties I: locally finite varieties*, Algebra Universalis **32** (1994), 224–269.
[2] D. Hobby and R. McKenzie, *The Structure of Finite Algebras*, Contemporary Mathematics v. 76, American Mathematical Society, 1988.
[3] K. A. Kearnes, *Type restriction in locally finite varieties with the* CEP, Canadian Journal of Mathematics **43** (1991), 748–769.
[4] K. A. Kearnes, *Cardinality bounds for subdirectly irreducible algebras*, Journal of Pure and Applied Algebra **112** (1996), 293–312.
[5] K. A. Kearnes, *Varieties with a difference term*, Journal of Algebra **177** (1995), 926–960.
[6] K. A. Kearnes, E. W. Kiss, and M. Valeriote. *A geometrical consequence of residual smallness.* Manuscript, 1996.
[7] K. A. Kearnes and R. McKenzie, *Residual smallness relativized to congruence types*, preprint, 1996.
[8] E. W. Kiss, *An easy way to minimal algebras*, to appear in the International Journal of Algebra and Computation.

(Keith A. Kearnes) DEPARTMENT OF MATHEMATICAL SCIENCES, UNIVERSITY OF ARKANSAS, FAYETTEVILLE, ARKANSAS 72701, USA

(Emil W. Kiss) EÖTVÖS UNIVERSITY, DEPARTMENT OF ALGEBRA AND NUMBER THEORY, 1088 BUDAPEST, MÚZEUM KRT. 6–8, HUNGARY
*E-mail address*, Keith A. Kearnes: `kearnes@comp.uark.edu`
*E-mail address*, Emil W. Kiss: `ewkiss@cs.elte.hu`